\pdfoutput=1
\documentclass[journal]{IEEEtran}
\usepackage{graphicx}
\usepackage[cmex10]{amsmath}
\usepackage{amsfonts}
\usepackage{amssymb}
\usepackage{xspace}
\usepackage{color}
\usepackage{subfig}
\usepackage{algorithm2e}
\usepackage{url}
\usepackage{hyperref}
\usepackage{etoolbox}

\graphicspath{{./images/}}

\newcommand{\Omit}[1]{}

\newif \ifindi\inditrue

\DeclareMathOperator*{\argmin}{arg\,min\:}

\allowdisplaybreaks 

\newtoggle{future}
\toggletrue{future}

\title{Dynamic Optimal Power Flow in Microgrids using the Alternating Direction Method of Multipliers}
\author{Paul~Scott
		and Sylvie~Thi\'ebaux
\thanks{This work has been submitted to the IEEE for possible publication. Copyright may be transferred without notice, after which this version may no longer be accessible.  The authors are from the College of Engineering and Computer Science at the Australian National University and the Optimisation Research Group at NICTA.}}%

\begin{document}

\maketitle

\begin{abstract}
Smart devices, storage and other distributed technologies have the potential to greatly improve the utilisation of network infrastructure and renewable generation.
Decentralised control of these technologies overcomes many scalability and privacy concerns, but in general still requires the underlying problem to be convex in order to guarantee convergence to a global optimum. 
Considering that AC power flows are non-convex in nature, and the operation of household devices often requires discrete decisions, there has been uncertainty surrounding the use of distributed methods in a realistic setting.  This paper extends prior work on the alternating direction method of multipliers (ADMM) for solving the dynamic optimal power flow (D-OPF) problem.  We utilise more realistic line and load models, and introduce a two-stage approach to managing discrete decisions and uncertainty.  Our experiments on a suburb-sized microgrid show that this approach provides near optimal results, in a time that is fast enough for receding horizon control.  This work brings distributed control of smart-grid technologies closer to reality.
\end{abstract}


\section{Introduction}
\IEEEPARstart{T}{he} role of electricity market operators is to supply low cost and reliable electricity to customers.  This typically involves solving unit commitment (UC) and/or optimal power flow (OPF) problems at regular intervals.  In the classical form of these problems, loads are assumed to be inflexible, so only large generators and network components need to be considered.  However, this assumption will no longer hold true when smart devices, distributed generation, storage, and electric vehicles (EVs) are adopted on mass.  The traditional centralised market approaches were not designed to either operate where every consumer is an active participant, or handle their unique time-coupled behaviours.

In order to manage this huge increase in problem size, several works \cite{kraning2014,gatsis2012,mohsenian-rad2010a} have adopted distributed solving techniques.  In addition to greatly parallelising the problem, these distributed algorithms help to preserve the privacy of participants.  These algorithms require the problem to be convex in order to guarantee convergence to a globally optimal solution. However, the behaviour of many loads at a household level are discrete in nature \cite{ramchurn2011}, and the equations that govern how power physically flows on the network are non-convex.

This paper considers the task of balancing power within a microgrid.  We formulate the problem as a dynamic optimal power flow (D-OPF) problem to account for multiple time steps over a day, and solve it in a distributed manner by adopting the alternating direction method of multipliers (ADMM) approach presented in \cite{kraning2014}.  We extend the method to more accurate power flow models and introduce a two-stage pricing mechanism to manage integer variables and uncertainty.  We find that this approach achieves near optimal solutions in a time frame that is fast enough to work with receding horizon control.  This brings the overall approach closer to the point where it can be deployed on a real network.

The introductory sections \ref{sec:related}-\ref{sec:microgrid} of this paper cover related work, our network model, the ADMM method and describe our test microgrid.  This is followed by three sections \ref{sec:powerflow}-\ref{sec:uncertainty} that focus on each of our contributions.

\section{Related Work}\label{sec:related}
Demand response (DR) is the name often given to the control of distributed technologies.  Much of the existing work on DR has focused on using real-time pricing (RTP) as a control signal \cite{mohsenian-rad2010b,chen2012,scott2013,tischer2011}.  In these methods, participants receive a RTP signal and individually optimise their own behaviour, so as to minimise a combination of monetary and discomfort costs.  Other approaches have utilised non-pricing control signals, which are simpler to implement, but are limited in the types of loads that they can model \cite{vandenbriel2013,shinwari2012}.

These approaches can be thought of as a form of open loop control, because the agent that sets the control signal (RTP or otherwise), at best can only estimate how consumers will respond to it.  The work in \cite{gatsis2012} takes the idea of RTP and closes the loop.  The RTP is iteratively updated by a central agent, with consumers communicating their best responses to the price before implementing them.  An alternative iterative procedure is introduced in \cite{mohsenian-rad2010a}, where consumers cooperate to reduce total generation costs in a distributed manner. This approach is robust to gaming, but individual households only have a small incentive to change their habits, as their costs are largely independent of their own behaviour.

The approaches discussed so far do not model the electricity network, so cannot account for real power losses, reactive power, voltage limits and line thermal limits.  Most of the work on distributed algorithms that do model the network have used ADMM as a solving technique, due to its ease in decomposition, and its convergence guarantees on a wide range of problems \cite{boyd2011}.
One of the first authors to apply ADMM to power networks was Kim et al.\ \cite{kim2000}. They decomposed an OPF problem into regions and compared ADMM to two other approaches.  They found that ADMM provided significant speed improvements over a centralised approach, and that privacy is preserved between regions.  In \cite{phan2014}, Phan et al.\ used ADMM to improve the parallel solve time for a security-constrained OPF (SCOPF) problem.  They decomposed across scenarios and found it to have comparable solve times to a Bender's cut approach in many instances.

The decomposition power of ADMM was taken one step further when Kraning et al.\ \cite{kraning2014} decomposed all network components for a D-OPF problem.  This is effectively a principled method for settling RTPs for each bus, also known as locational marginal prices.  Their experiments showed that very large problems could be solved efficiently in a parallel environment. In this paper we extend this work to include more sophisticated line models, and a method for managing discrete variables and uncertainty, thereby making ADMM a practical approach to modern D-OPF problems.

\section{Network Model}
The overall objective is to minimise the cost of supplying electricity.  We formulate this as D-OPF problems embedded within a receding horizon control process, in order to accurately control time-coupled components in an uncertain environment.  In receding horizon control, a D-OPF is first solved over a horizon of $n \in \mathbb{N}$ time steps, the decision in the first step is acted on, and then the process repeats with the window shifted forward by one. In this paper we focus on solving the D-OPF within a single horizon, and the actions that agents take to implement the first decision.  We formalise the network model with this problem in mind.

\begin{figure}[t]
\centering
\includegraphics[width=\linewidth]{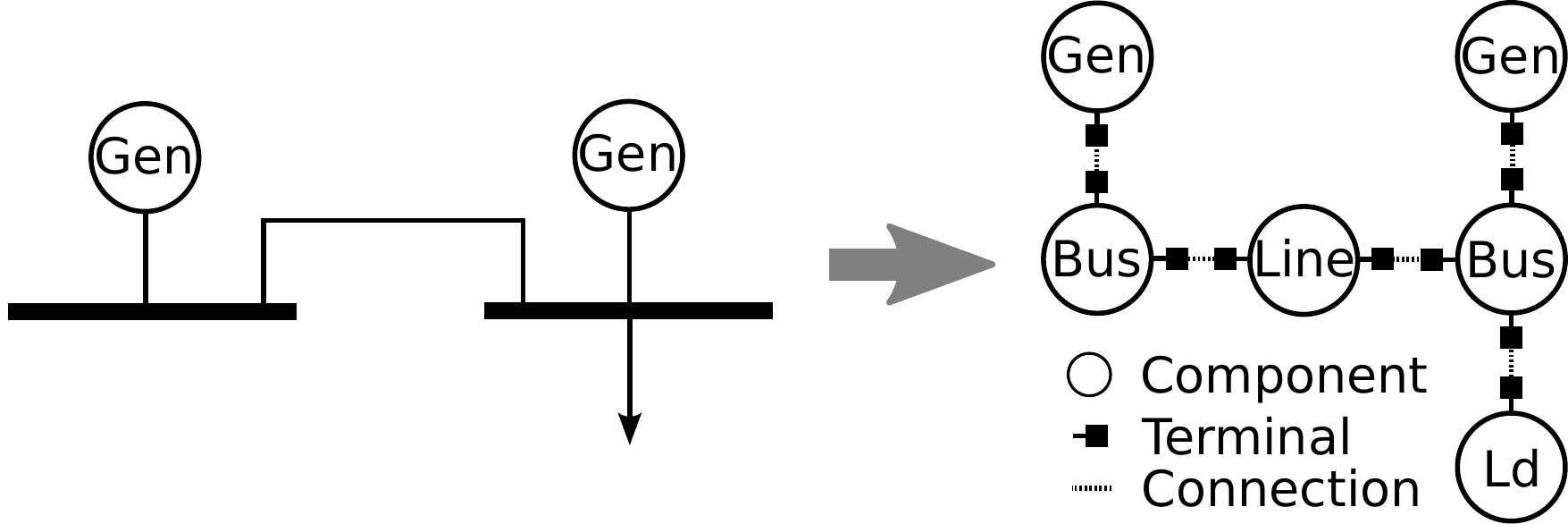}
\caption{Conversion to component orientated representation.}
\label{fig:network_desc}
\end{figure}

Note that the notation and formulation that we use is different from what is standard in power systems.  This is necessary in order to decompose the network and distribute the problem.  Fig.~\ref{fig:network_desc} highlights the difference between a standard line diagram and our formulation.

A network $N$ consists of a set of components $C$, terminals $T$ and connections $L$.  Each component $c \in C$ (e.g., bus, line, generator, load) has a set of terminals $T_c \subseteq T$ which can be connected to the terminals of other components, where the $T_c$ sets partition $T$.
Each connection $l\in L$ is a pair of terminals, i.e.\ $L \subseteq T\times T$.

\subsection{Connections}
A terminal represents a point of connection for an electrical component.  We use the quantities of real power, reactive power, voltage and voltage phase angle ($p, q, v, \theta \in \mathbb{R}^n$ respectively) to model the flow of power into a component through a terminal.  These are vectors in order to capture each time step in the horizon.  For convenience we use a parent vector $y_i \in \mathbb{R}^{4n}$ to represent all variables for a terminal $i \in T$, where $y_i := (p_i, q_i, v_i, \theta_i)^\mathsf{T}$.  When two terminals are connected together, $(i, j) \in L$, we pose the following constraints:
\begin{align}
p_i + p_j = 0,\quad q_i + q_j = 0,\quad v_i - v_j = 0,\quad \theta_i - \theta_j = 0\nonumber
\end{align}
The first two constraints ensure the power that leaves one terminal enters the terminal it connects to.  The second two constraints ensure that the terminals have the same voltage and phase angle.  This duplication of variables is necessary in order to decompose the problem for the ADMM algorithm.  To avoid confusion, recall that connections and terminals are different from lines and buses (see Fig.~\ref{fig:network_desc}).

We rewrite these constraints as $y_i + A y_j = 0$ for $y$, where $A$ is the appropriate $4n\times4n$ diagonal matrix.  Further, we define the function $h: \mathbb{R}^{4n}\times \mathbb{R}^{4n} \mapsto \mathbb{R}^{4n}$ as the LHS of this constraint for convenience: $h(y_i,y_j) := y_i + A y_j$.

\subsection{Components}
At a high level, each component $c \in C$ has a variable vector $x_c \in \mathbb{R}^{a_c}$, an objective function $f_c: \mathbb{R}^{a_c} \mapsto \mathbb{R}$, and a constraint function $g_c: \mathbb{R}^{a_c} \mapsto \mathbb{R}^{b_c}$, where $g_c(x_c)\leq 0$.  For a component $c\in C$, the vector $x_c$ includes all terminal variables for that component: $y_i, \forall i \in T_c$.

In the following sections we describe at a lower level the models for a variety of components.
When necessary we use $t \in \{0,\ldots,n\}$ to index vectors by time, otherwise we imply standard vector operations. The index where $t=0$ is used to represent the value of the variable at the beginning of the current horizon, which we assume is known.

\subsubsection{Bus}
A bus has a variable number of terminals which depends on how many other components connect to it. For example, a bus might be connected to a generator, a load and 3 lines for a total of 5 terminals.  Regardless of the number of terminals, the constraints take the form:
\begin{align}
&\sum_{i \in T_c} p_i = 0\quad \sum_{i \in T_c} q_i = 0\nonumber\\
&\forall i, j \in T_c: v_i = v_j,~\theta_i = \theta_j\nonumber
\end{align}
The first two constraints are an expression of Kirchhoff's current law (KCL) in terms of power flows.  The remaining constraints ensure that all terminal voltages and phase angles are the same.

\subsubsection{Line}\label{sec:line}
A line is a two terminal component which transports power from one location to another, typically from bus to bus.  We model a line as having a constant conductance $g\in\mathbb{R}_+$, susceptance $b\in\mathbb{R}$ and maximum apparent power $s\in\mathbb{R}_+$. The AC power flow equations are derived from Ohm's law, where $\forall i,j\in T_c,i\neq j$:
\begin{align}
&p_i = g v_i^2 - g v_i v_j \cos(\theta_i - \theta_j) - b v_i v_j \sin(\theta_i - \theta_j)\label{eq:ac_beg}\\
&q_i = -b v_i^2 + b v_i v_j \cos(\theta_i - \theta_j) - g v_i v_j \sin(\theta_i - \theta_j)\\
&s^2 \geq p_i^2 + q_i^2,\quad\underline{v} \leq v_i \leq \bar{v},\quad\theta_i - \theta_j\leq\bar{\theta}\label{eq:ac_end}
\end{align}
These constraints are identical for each time step, so we have left out the indexing by time to improve clarity.
These equations are non-convex, so they are often either approximated or relaxed, as we will discuss further in Section~\ref{sec:powerflow}.

\subsubsection{Generator}
A generator is a single terminal component which produces real and reactive power.  In our formulation the generator has a floating phase angle and voltage.  A generator has lower and upper real and reactive power limits such that $p_{i,t} \in [\underline{p},\bar{p}]$ and $q_{i,t} \in [\underline{q},\bar{q}]$, a ramping rate $p^r\in\mathbb{R}_+$ and a quadratic cost function $f$ for generation costs:
\begin{align}
&f(x) = p_i^\mathsf{T}\Psi p_i - \psi^\mathsf{T} p_i\nonumber\\
&\forall t\in\{1,\ldots,n\}:-p^r \leq p_{i,t} - p_{i,t-1}  \leq p^r\nonumber
\end{align}
where $\Psi\in\mathbb{R}_+^{n\times n}$ is a diagonal matrix and $\psi\in\mathbb{R}_+^n$.

\subsubsection{Shiftable Load}\label{sec:shiftable}
A shiftable load is a single terminal component used to model electrical loads like dish washers and clothes dryers.  These loads must start running between an earliest and a latest start time: $t^{e},t^{l}\in\mathbb{N}$.  To model this we introduce binary variables $u \in \{0,1\}^n$ for the horizon.  A value of $1$ indicates that the component starts at the given time.  A component runs for a duration of $d\in\mathbb{N}$ consecutive time steps, during which it consumes a load of $p^\mathrm{nom}\in\mathbb{R}$.
\begin{align}
&p_{i,t} = p^\mathrm{nom}\sum_{t' = t-d+1}^tu_{t'}\quad \sum_{t=t^e}^{t^l} u_t = 1\nonumber\\
&\forall t \not\in \{t^e, \ldots, t^l\} : u_t = 0\nonumber
\end{align}
A convex relaxation of this component can be obtained by relaxing the integrality requirement: $u \in [0,1]^n$.

\subsubsection{Battery}
A battery is a single terminal component with stored energy $E\in[0,\bar{E}]^n$.  The charge $p^c\in \mathbb{R}_+^n$ and discharge $p^d\in \mathbb{R}_+^n$ variables transfer energy to/from the battery:
\begin{align}
&E_t = E_{t-1} + \eta p_t^c - p_t^d,\quad E_n \geq \bar{E}/2\nonumber
\end{align}
where the constant $\eta\in[0,1]$ gives the efficiency of the battery and the net power into the component is given by $p_i = p^c - p^d$.  The second constraint above prevents a greedy discharge of power at the end of the horizon.


\section{Alternating Direction Method of Multipliers}\label{sec:admm}
The goal is to minimise the sum of all component cost functions, subject to component and terminal connection constraints.  This is a utilitarian view of the problem.
\begin{align}
&\min_x\sum_{c \in C} f_c(x_c)\nonumber\\
&\text{ s.t. } \forall c \in C: g_c(x_c) \leq 0 \nonumber\\
&\phantom{\text{ s.t. }} \forall (i, j) \in L: h(y_i, y_j) = 0\nonumber
\end{align}
The augmented Lagrangian applied to the connection constraints is given by:
\begin{align}
\mathcal{L}(L, y, z, \lambda, \rho) &:= \sum_{c \in C} f_c(x_c)\nonumber\\
	&+ \sum_{(i, j) \in L}
	\left[\frac{\rho}{2} \|h(y_i, z_j)\|_2^2 +
	\lambda_{i,j}^\mathsf{T} h(y_i, z_j)
	\right]\nonumber
\end{align}
where $\rho\in(0,\infty)$ is a penalty parameter and $\lambda_{i,j} \in \mathbb{R}^{4n}$ are the dual variables for the connection constraints.  The dual variables also represent the locational marginal prices in our problem.


The alternating direction method of multipliers (ADMM) is a variation of the standard augmented Lagrangian method that enables problem decomposition \cite{boyd2011,douglas1956,gabay1976}. A single iteration consists of two phases followed by a dual variable update.  Components are each allocated to one of the two phases.  The component sets $C_1$ and $C_2$, and the variable vectors $x_1$ and $x_2$ represent this allocation.

The connections are split into three parts: $L_1$, $L_2$ and $L_{1,2}$. The intra-phase connections $L_1$ ($L_2$) are those that are between components in $C_1$ ($C_2$).  The inter-phase connections $L_{1,2}$ are those where one component is in $C_1$ and the other is in $C_2$.

The superscript $k\in\mathbb{N}$ is used to indicate the $k$-th iteration.  At the start of the algorithm all terminal and dual variables are initialised to some values $y_i^{(0)}$ and $\lambda_{i,j}^{(0)}$.  For the $k$-th iteration ADMM proceeds as follows:
\begin{enumerate}
\item Optimise for $x_1$, holding $x_2$ constant at its $k-1$ value
\item Optimise for $x_2$, holding $x_1$ constant at its $k$ value
\item Update the dual variables $\lambda$
\end{enumerate}
For our optimisation problem this becomes:
\begin{align}
&x_1^{(k)} = \argmin_{x_1} \mathcal{L}(L_{1,2}, y, y^{(k-1)}, \lambda^{(k-1)}, \rho^k)\label{eq:admm_st}\\
&\quad\text{ s.t. } \forall c \in C_1: g_c(x_c) \leq 0 \nonumber\\
&\quad\phantom{\text{ s.t. }} \forall (i, j) \in L_1: h(y_i, y_j) = 0\nonumber\\
&x_2^{(k)} = \argmin_{x_2} \mathcal{L}(L_{1,2}, y^{(k)}, y, \lambda^{(k-1)}, \rho^k)\\
&\quad\text{ s.t. } \forall c \in C_2: g_c(x_c) \leq 0 \nonumber\\
&\quad\phantom{\text{ s.t. }} \forall (i, j) \in L_2: h(y_i, y_j) = 0\nonumber\\
&\forall (i, j) \in L_{1,2}: \lambda_{i,j}^{(k)} = \lambda_{i,j}^{(k-1)} + \rho^{(k)} h(y_i^{(k)}, y_j^{(k)})\label{eq:admm_en}
\end{align}

In the simple case when $\rho$ is constant, $f_c$ and $g_c$ are convex, and $h$ is affine, ADMM converges to a global optimum \cite{boyd2011}.

If a component has no intra-phase connections, then its subproblem can be solved independently of all other components in the same phase.  We adopt the partitioning scheme where $C_2$ contains all buses and $C_1$ the rest of the network.  This allows us to fully separate all components, since buses will never connect to other buses ($L_2 = \emptyset$) and non-bus components will never connect to other non-bus components ($L_1 = \emptyset$).  As an additional benefit, some components are simple enough that we can analytically calculate the solution to their subproblem at each iteration, instead of invoking an optimisation routine.  We adopt such an analytical approach for buses as proposed in \cite{kraning2014}.  Other partitioning schemes are discussed in Section~\ref{sec:powerflow}.

\subsection{Residuals and Stopping Criteria}
As in \cite{kraning2014}, we use primal and dual residuals to define the stopping criteria for our algorithm.  The primal residuals represent the constraint violations at the current solution.  We combine the residuals of all connections into a single vector $r_p$.  By indexing into the inter-phase connections $L_{1,2} = \{(i_1, j_1), (i_2, j_2), \ldots \}$, the primal residual is given by:
\begin{equation}
r_p^{(k)} := (h(y_{i_1}^{(k)}, y_{j_1}^{(k)}), h(y_{i_2}^{(k)}, y_{j_2}^{(k)}), \ldots)^\mathsf{T}\nonumber
\end{equation}
The dual residuals give the violation of the KKT stationarity constraint at the current solution.  We collect the dual residuals for each each connection into the vector $r_d$.  For ADMM, the dual residuals are given by (see \cite{boyd2011} for derivation):
\begin{equation}
r_d^{(k)} := (\lambda_{i_1,j_1}^{(k)} - \lambda_{i_1,j_1}^{(k-1)}, \lambda_{i_2,j_2}^{(k)} - \lambda_{i_2,j_2}^{(k-1)}, \ldots)^\mathsf{T}\nonumber
\end{equation}
These residuals approach zero as the algorithm converges to a KKT point.  We consider that the algorithm has converged when the scaled $2$-norms of these residuals are smaller than a tolerance $\epsilon$: $\frac{1}{\sqrt{M}}\|r_p^{(k)}\|_2 < \epsilon$, $\frac{1}{\sqrt{M}}\|r_d^{(k)}\|_2 < \epsilon$.  Here $M$ is the total number of inter-phase terminal constraints $4n|L_{1,2}|$ minus the number of terminal constraints that are trivially satisfied (e.g., floating voltages and phase angles for generators).  It is used to keep the tolerance independent of problem size.


\section{Implementation}
We developed an experimental implementation of the ADMM method in C++ using Gurobi \cite{gurobi2014} and Ipopt \cite{wachter2006,hsl2014} as backend solvers.  Gurobi is used for mixed-integer linear or quadratically constrained problems, and Ipopt for more general nonlinear problems.  CasADi \cite{andersson2013} was used as a modelling and automatic differentiation front end to Ipopt.
This implementation was designed with flexibility in mind, so that a wide range of experiments could be conducted.  For a more specialised and performance orientated implementation see \cite{kraning2014}.
Although this is a sequential implementation of the ADMM method, we timed the slowest component at each iteration to get an idea of how long a fully distributed implementation would take.

The experiments were run on machines with 2 AMD 6-Core Opteron 4184, 2.8GHz, 3M L2/6M L3 Cache CPUs and 64GB of memory.

\section{Test Microgrid}\label{sec:microgrid}
The distributed algorithm we have presented can be used to control all levels of the electricity system, potentially in a hierarchical manner \cite{kraning2014}. We focus our attention on a microgrid-like distribution network since that is where we expect to have the greatest impact with the adoption of distributed generation and storage.

Our experiments are based around a modified 70 bus 11kV benchmark distribution network \cite{das2006}, which was chosen due to its comparable size to a suburb.  We close all tie lines in the network in order to change it from a radial to a meshed configuration.  This is because we expect microgrids to take on more of a meshed network structure to improve reliability, efficiency and to better utilise distributed generation.

The benchmark comes with a static PQ load at each bus, which we replace with a number of houses, depending on the size of the load.  The houses are connected directly to the 11kV buses as we have no data on the low voltage part of the network.  We assume that the power bounds we place on each household will be sufficient to prevent any capacity violation of the low voltage network.

A house is an independent agent that manages subcomponents. These include an uncontrollable background power draw, two shiftable loads, and optional battery and solar PV systems.  A house has a single terminal through which it can exchange real and reactive power with the rest of the network.  Each house has an apparent power limit of $s=10$kVA.

We develop a typical house load profile $l_t$ by modifying an aggregate Autumn load profile for the ACT region in Australia (data from \cite{aemo}).  We assume that households consume on average 20kWh per day.  This provides the basis for all uncontrollable household background loads.  For the purposes of these experiments, we assume that the static PQ loads in the benchmark were recorded when loads were at 75\% of their peak.  We divide the benchmark static real power at each bus by how much power a typical house consumes at 75\% of its peak power (1.45kW).  Rounding down this number gives us an estimate of the number of houses which would be located at a given bus.  This approach produces a total of $h=3674$ houses for the network, about the size of an Australian suburb.

We place two generators in the network where the distribution system connects to upstream substations.  These can be thought of as either dispatchable microgenerators or as representing the cost of importing power into the microgrid.

We randomise some of the generator and household load parameters to produce different problem instances, as can be seen in Table~\ref{tab:micro_pars}. The time horizon spans 24 hours with 15 minute time steps, which produces a problem instance with over 2 million variables per horizon. The experiments were run with a primal and dual stopping tolerance of $\epsilon=10^{-4}$ and a fixed penalty parameter of $\rho=0.5$.  To improve numerical stability, we scale the system to a per-unit representation with base values at 11kV and 100kVA.  This means that a real power residual of $10^{-4}$ translates to 10W for a connection, or about 1\% of the average household load.

\begin{table}[t]
\centering
\caption{Component parameters.}
\label{tab:micro_pars}
\begin{tabular}{r|r|r|r}
Component & Param & Value & Units\\
\hline
Generator& $\psi_t$ & $\sim\max(4, \mathcal{N}(40, 8^2))$ & \$/MWh\\
Generator& $\Psi_{t,t}$ & $\sim\max(1, \mathcal{N}(10, 2^2))$ & \$/MWhMW\\
Generator& $\underline{p},\bar{p}$ & $-s\times h /2,0$ & kW\\
Generator& $\underline{q},\bar{q}$ & $0.2\underline{p}, -0.2\underline{p}$ & kvar\\
House& $p_t$ & $\sim\mathcal{N}(l_t, (0.2l_t)^2)$ & kW\\
House& $q_t$ & $0.3p_t$ & kvar\\
Shift 1& $d$ & $\sim\max(15, \mathcal{N}(90, 18^2))$ & min\\
Shift 1& $p^\mathrm{nom}$ & $\sim\max(0.3, \mathcal{N}(3, 0.6^2))$ & kW\\
Shift 2& $d$ & $\sim\max(15, \mathcal{N}(60, 12^2))$ & min\\
Shift 2& $p^\mathrm{nom}$ & $\sim\max(0.1, \mathcal{N}(1, 0.2^2))$ & kW\\
Shift& $t^e,t^l$ & $0, n-d$ &
\end{tabular}
\end{table}

In addition to the 70 bus network described here, we also ran a series of experiments on randomly generated networks similar to those described in \cite{kraning2014}. These randomly generated networks ranged in size from 20 to 2000 buses, and were designed to be highly congested.  We will occasionally mention some of the results from these random networks when they differ from our 70 bus microgrid.

\section{Impact of Power Flow Models}\label{sec:powerflow}
In this section we investigate how the ADMM method performs with different power flow models.  We assess 4 alternative models to the one presented in \cite{kraning2014}, all of which give more accurate results, but are slower to converge.

\subsection{Power Flow Models}
Due to their non-convex nature, the AC power flow equations (\ref{eq:ac_beg}--\ref{eq:ac_end}) are often either relaxed or approximated.  Convex relaxations include a quadratic constraint (QC) model \cite{hijazi2014}, a semi-definite program \cite{bai2008}, the dist-flow (DF) relaxation \cite{farivar2011,baran1989c} and an equivalent SOCP relaxation \cite{jabr2006,bose2014}.  Approximations include the linear DC (DC) model \cite{schweppe1970,stott2009}, the LPAC model \cite{coffrin2014} and the equations proposed by Kraning et al.\ (K) \cite{kraning2014}.  These alternative models are often used to solve difficult power network optimisation problems, for example, OPF, OPF with line switching, capacitor placement and expansion planing.  ADMM was used with the DF model in \cite{sulc2014} in order to control the reactive power of inverters on a radial network, and the K model was used with ADMM in \cite{kraning2014} to solve the D-OPF problem.

What is lacking is an understanding of the relative strengths and weaknesses of line models when used in a distributed algorithm.  In this section we compare how the distributed ADMM algorithm performs when using the AC, QC, DF, DC, and K line models.  We compare the differences in solution quality, feasibility, processing time and number of iterations for our test network.

\subsection{Experiments}
We generated 60 instances of our test network with the binary variables for the shiftable devices relaxed.  To evaluate the quality of the solutions we calculate the percentage difference in objective value relative to the best known AC solution: $100\cdot(f - f_{best})/f_{best}$.  The means and standard deviations of the 60 instances are:
\begin{center}
\begin{tabular}{r r|r r}
QC & -0.031\%~(0.008\%) & DF & 0.039\%~(0.018\%)\\
DC & -3.541\%~(0.072\%) & K  & 4.726\%~(0.090\%)\\
\end{tabular}
\end{center}
We don't have a guarantee that the AC solution is optimal because the equations are non-convex. However, the relaxed models (QC, DF) provide a lower bound on the global optimum, allowing us to identify in the worst case how far the AC solution is from optimality.  The DF model produces results that are slightly above the AC model, when it should be providing a lower bound.  This is because the result is within the margin of error for our stopping criteria, which we estimated to be 1\%.  By decreasing $\epsilon$ we found the positive difference to shrink further into insignificance.  Based on our results the optimality gap is tiny for both the QC and DF models.  These results give us confidence that the non-convex AC model, which is the only one that guarantees Ohm's law is satisfied, produces solutions that are very close to optimal. Other work has come to a similar conclusion in a different setting \cite{hijazi2014,phan2014}. 

The DC model underestimates the optimal value by around 3.5\% while the K model overestimates it by around 4.7\%.
Part of the reason for this is that the DC model ignores line losses and the K model is overestimating line losses.


Table~\ref{tab:lines_time} provides the number of iterations and time taken to converge in the form of means and standard deviations.  The parallel solve time is the amount of time required to solve the problem in a fully distributed implementation.  This was measured by summing together the time of the slowest component at each iteration.  In absolute terms the parallel solve time is relatively small despite the fact that our implementation was designed with flexibility in mind, not performance.  That said the K model is significantly faster relative to the other models.  It converges in half the number of iterations required by the next fastest model, but as we have seen it gives us an inaccurate result.

\begin{table}[t]
\centering
\caption{Iterations and parallel solve time for line models.}
\label{tab:lines_time}
\begin{tabular}{r r|r}
& Iterations (std.) & Time in sec (std.)\\
\hline
AC &  1945~(17) & 148~(12)\\
QC &  1951~(14) & 546~(33)\\
DF &  1933~(26) & 110~(8)\\
DC &  4140~(50) & 244~(8)\\
K  &  1027~(52) & 15~(1)\\
\end{tabular}
\end{table}


The highly congested random networks that we generated produce similar results.  However, for a number of instances the K model was infeasible (would not converge) where we had a valid AC solution, due to the exaggeration of line losses.

\begin{figure}[t]
\centering
\includegraphics[width=0.8\linewidth]{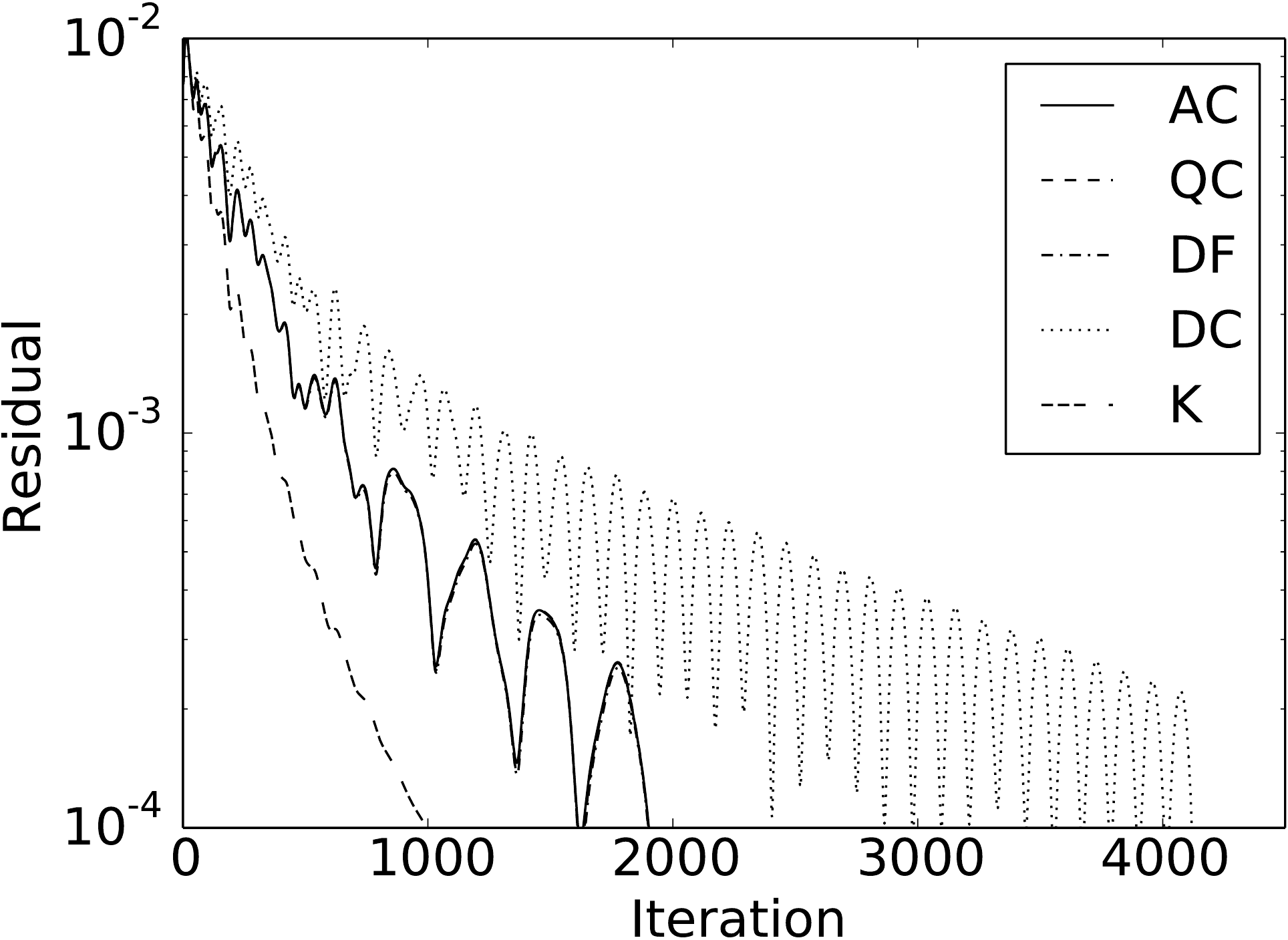}
\caption{Convergence of primal residuals.}
\label{fig:lines_res}
\end{figure}

Fig.~\ref{fig:lines_res} shows an example of the primal residual convergence for different line models (the dual residuals are similar).  The AC, QC and DF models overlap.  One unintuitive result is the fact that the DC model converges poorly when it is in fact a very simple linear model.  Large oscillations build up during the solution process which slows the rate of convergence.  We performed a series of experiments in order to get a better understanding of this effect. The best explanation we have is that the DC model behaves like an undamped system, as it has no line losses and only linear constraints.  The DC model will respond stronger for a given change in its terminal dual variables.  The net effect is that oscillations build up across the network during the solving process.  On the other hand the K model overestimates line losses, which means it is much less sensitive and no oscillations form.  The AC, QC and DF models are somewhere in between these two extremes.

In reality communication delay will play a major part in the total solve time for the algorithm.  For example, if we assume that it takes 60ms to communicate messages on our network (which happens a minimum of two times per iteration), then 1000 iterations would require 2 minutes of communication time.  In certain circumstances it may be beneficial to cut down the total number of iterations, even if it requires more processing time per iteration.

We experimented with clustering some lines and buses together instead of fully decomposing the network.  The idea is that, depending on the communication delay, this could reduce overall solve time by solving more of the problem at once.  This was clearly beneficial for the random networks, but for the 70 bus microgrid the results so far are inconclusive.  For example, clustering the 70 bus network into 16 parts halved the number of iterations required, but the overall solve time worsened.  We expect this to change once we start to optimise our implementation, by reducing solver overhead and decomposing network components across time steps.

It is important to point out that we are giving the algorithm a naive starting point for both the primal and dual variables.  In practice, the receding horizon control scheme will provide an excellent warm starting point, because the values from the previous horizon can be used for all but one time step.  As a sanity check, we performed warm starting experiments for the AC model.  Similar to what was done in \cite{kraning2014}, we duplicate a problem instance and then randomly resample the household background power and shiftable device power parameters according to the rule: $p \sim p\mathcal{N}(1,\sigma^2)$.  We used the solution of the original instance as a starting point for the modified instance.  For $\sigma=0.2$ the warm started run only needed 11\% of the original iterations on average.  In a second experiment we fully correlated the resampling step, which could represent a correlated change in solar panel output for many households.  With $\sigma=0.2$, only 29\% of the original iterations were required on average.

These results show the feasibility of using the non-convex AC power flow equations for solving a distributed D-OPF problem.  The K model converges much faster, but it is unlikely to be usable in a realistic setting, as it ignores voltages and reactive power, and produces overly high costs.  For these reasons, we adopt the AC model for the rest of our experiments.

\iftoggle{future}{}{%
In future work we will investigate if a reduction in total solve time can be gained by running the K model for a certain duration and then switching to the AC model.  Also, there is the potential to further parallelise the problem by decomposing across time within certain components.  We will also consider further improving convergence rates by investigating alternatives to the ADMM method, such as those that use second-order information.
}


\section{Discrete Decisions}
We now extend the algorithm so that it can deal with discrete decisions.  Components with discrete decisions are be common in realistic household models and they are required to model more complicated generators and network switching events.  We develop a framework that can accommodate discrete variables otherwise we could risk network violations through unpredictable household behaviour.

\subsection{Methods}
We investigate three tractable methods for dealing with integer variables which have no global optimality guarantees.  Just as we did for the AC equations, we will compare our result to a lower bound in order to get an understanding of the optimality gap.  We categorise these methods as:
\begin{itemize}
\item Relax and price (RP)
\item Relax and decide (RD)
\item Unrelaxed (UR)
\end{itemize}
The RP and RD approaches are broken up into two stages. The first stage, also known as the negotiation stage, is common to both methods.  Here all integer variables are relaxed and the distributed algorithm is run until convergence.  At this point the integer variables may take on fractional values, and this solution gives a lower bound on the optimal.  In the second stage each component makes a local decision in order to force any fractional values to integers.

Recall from Section~\ref{sec:shiftable} that shiftable devices have a binary variable $u_t$ for each time step, only one of which can take on the value 1 to indicate the starting time.  In the second stage of the RP method, each house performs a local optimisation to determine how to enforce integer feasibility of $u_t$.  We designed a range of cost functions which penalise a component if it changes its terminal values from those negotiated.  For a given cost function each house solves a Mixed-Integer Program (MIP) to obtain an integer-feasible solution.  The two most effective cost functions that we identified are:
\begin{align}
f_0(y, \hat{y}, \hat{\lambda}) &=  \hat{\lambda}^\mathsf{T} y + \alpha h(y, \hat{y})^\mathsf{T}h(y, \hat{y})\\
f_3(y, \hat{y}, \hat{\lambda}) &= \hat{\lambda}^\mathsf{T} A\hat{y} + \alpha h(y, \hat{y})^\mathsf{T}\Lambda h(y, \hat{y})
\end{align}
Where for a given house to bus connection $\hat{y}$ is the negotiated terminal values for the bus and $\hat{\lambda}$ the negotiated dual variables.  We use $\Lambda$ to represent the diagonal matrix where $\Lambda_{i,i} := |\hat{\lambda}_i|$ and $\alpha$ is a penalty parameter.

The first function charges households at the negotiated price for what they \emph{actually} consume, but they are also charged a quadratic penalty for operating away from the negotiated consumption.  The second function requires the household to pay for all power that was negotiated in the first stage.  Like the first function a penalty is charged for operating away from the negotiated operating point, however the penalty is scaled by the dual variables which can vary with time.

After this local optimisation step, we check that the solution is feasible and what the overall cost is. In order to do this we need to put some degrees of freedom back into the problem. In power networks the dispatch of generators are established in advance in response to an estimated demand. This forecast is never perfect, so a certain number of generators are paid to perform frequency regulation in order to balance demand in real time. In our experiments we employ both our generators for this use by allowing them to adjust their output. For these experiments we assume the same cost function and prices for both dispatch and frequency regulation.

In the second stage of the RD method, the largest $u_t$ value of each shiftable component is chosen to be fixed at 1 and the rest set to 0.  After this is done the distributed algorithm is started again in order to converge to a new solution that is integer feasible.

The final approach, UR, attempts to enforce integrality satisfaction at each iteration of the distributed algorithm.  We have already foregone theoretical convergence guarantees by our adoption of the non-convex AC equations.  Here we push the ADMM algorithm even further by allowing discrete variables into the algorithm (\ref{eq:admm_st}--\ref{eq:admm_en}), where Gurobi solves MIPs for houses, and Ipopt NLPs for lines.

\subsection{Experiments}
\begin{figure}[t]
\centering
\includegraphics[width=0.7\linewidth]{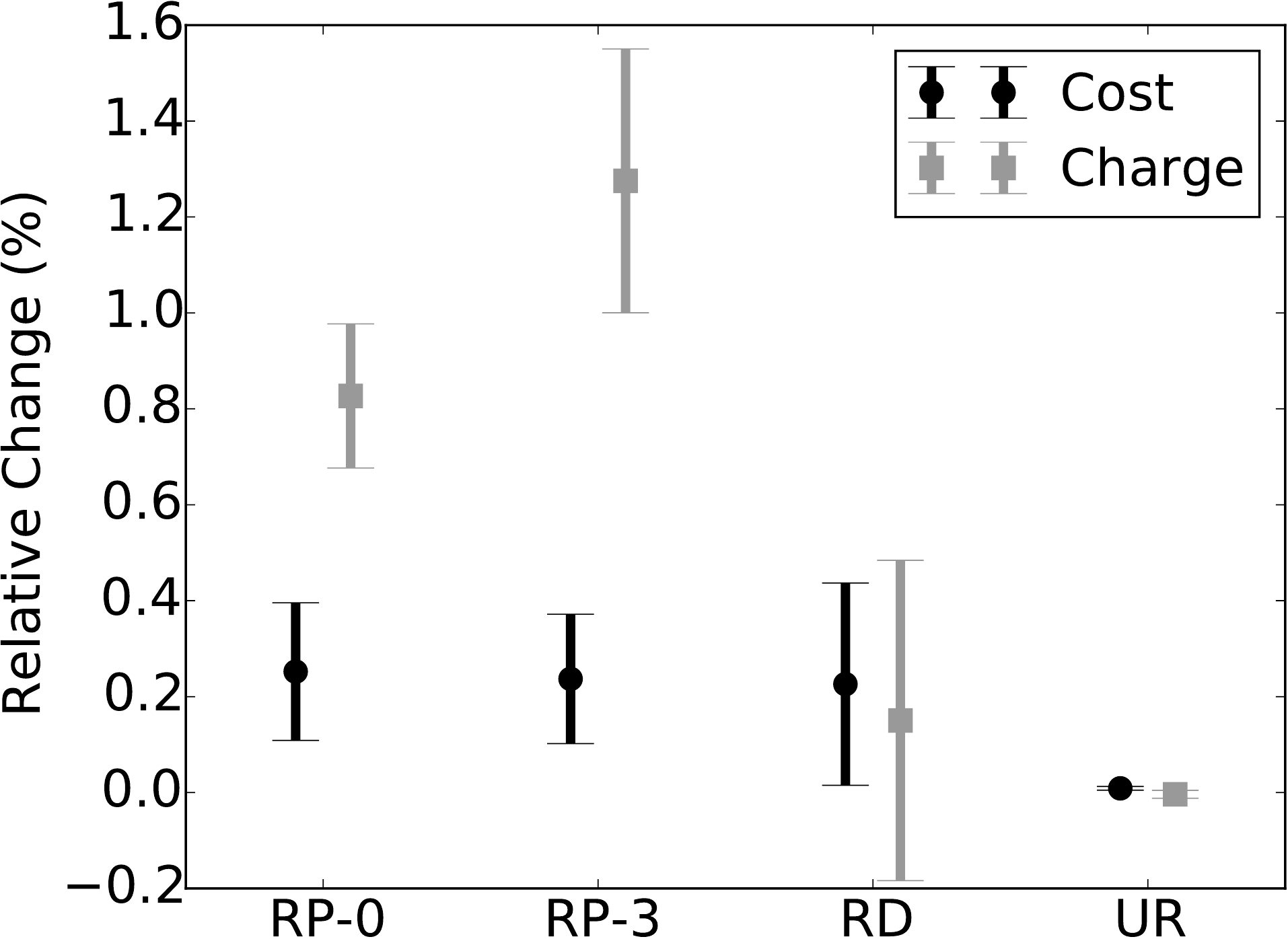}
\caption{Cost and charge error relative to relaxed solution.}
\label{fig:recost_comp}
\end{figure}

We ran experiments on 60 random instances of our test network with a penalty of $\alpha=10$. The results can be seen in Fig.~\ref{fig:recost_comp} which gives the fractional change with respect to the relaxed problem.  This is shown both in terms of the cost of generation and the charge to households.  The charge is the sum of household objective functions, which represents the amount of money they pay for their electricity.  For the RP methods this is given by the cost functions in the previous section.  For the RD and UR methods the charge is simply the final $\lambda^\mathsf{T}y$ for each house.

For the relaxed problem itself, the true cost of generation can be different from the amount households are charged, as we are dealing with marginal prices.  In addition to this, network congestion typically generates additional revenue above the cost of generation itself.  An increase in cost for the integer feasible solution relative to the relaxed problem is an indication of the additional cost to the generators for balancing supply.  Where household charge has increased relative to the relaxed solution, then households have decided to take on additional costs in order to change their consumption.

All methods produce costs that are within 1\% of the relaxed problem.  There is no significant difference between the four methods as they reside within our estimated margin of error based on our stopping tolerance.  What these results suggest is that we have a tight relaxation of the integer problem.  From other experiments, the relaxation of the shiftable devices is tight except where we have a heavily congested system.  This combined with the fact that each shiftable load only contributes a tiny amount to the overall power demand gives us a tight relaxation.

The RP method only marginally increases solve time above the results in the Section~\ref{sec:powerflow}. The RD method requires some of extra time as it performs a warm restart of the distributed algorithm.  The UR method took 1.7 times longer on average as it solves MIPs during each iteration.

The charges to households are significantly higher for the RP method without gaining any benefit in terms of reduced costs.  We ran the same experiments with a much smaller $\alpha$ which all but eliminated charges, without any increase to costs.  This suggests that for the sole purpose of managing discrete decisions, there is no need to have a strong penalty.

All of the methods we have presented provide an efficient means for dealing with the discrete decisions in a household. However, we will come to see in the next section why we favour the RP approach.

\iftoggle{future}{}{%
If we had larger discrete decisions, for example, shutting off large industrial plant, large generator start-up costs, or line switching, then we would likely have more problems with a continuous relaxation.  Future work will focus around investigating how the different methods handle more invasive discrete events and under different levels of network congestion.
}

\section{Pricing Uncertainty}\label{sec:uncertainty}
In this section we investigate the inclusion of stochastic components into our system.  Many parameters such as background household power consumption, solar PV output can only be estimated.  This means that the solutions that are negotiated through our distributed mechanism may not be applicable when it comes time to act.
For example, if the output of household solar PV systems is lower than expected, then certain lines in the network could be overloaded if the network was already running near capacity.

Part of the way of dealing with this is to run receding horizon control so that we only ever act on the most immediate time step before reoptimising. However even within the most immediate time step we still have to deal with uncertainty.  In this section we use the RP method from the previous section to encourage households to correct for any local sources of uncertainty.  Of the three methods designed to handle discrete variables, this is the only one that can provide this type of incentive.  To simplify the experiments we don't perform full receding horizon control.  With receding horizon control we expect the same trends, but with further improvement to costs.

\subsection{Experiments}
\begin{figure}[t]
\centering
\includegraphics[width=0.7\linewidth]{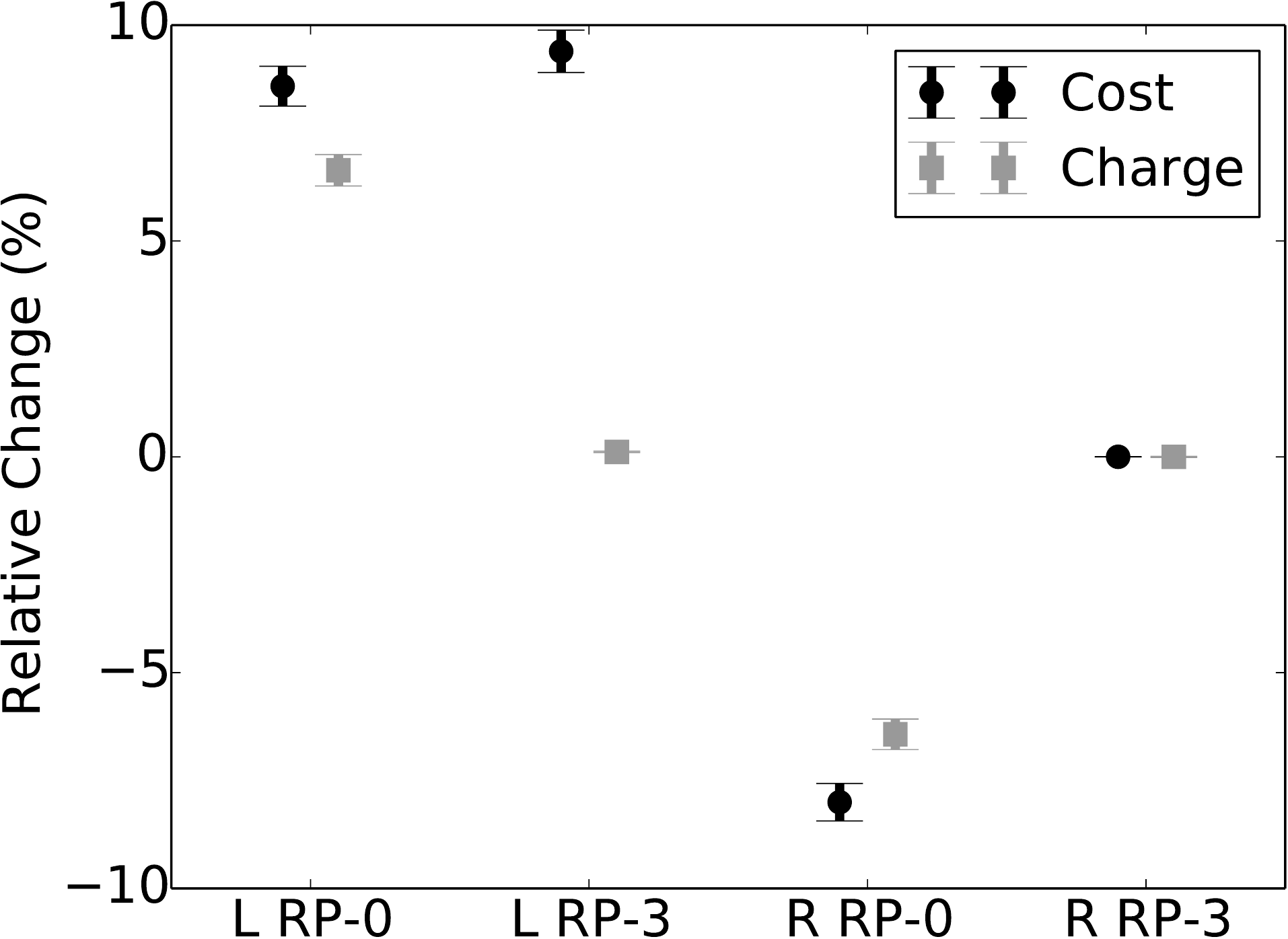}
\caption{Performance of cost functions with lowering (L) and raising (R) solar output.}
\label{fig:uncer_comp}
\end{figure}

We perform a simple experiment on the 70 bus microgrid where we have added 2kW PV systems and 2kWh batteries independently and at random to half of the houses.  The battery efficiencies $\eta$ were uniformly sampled from the interval $[0.85,0.95]$.  For normalised solar irradiance we use the simple relation: $I_t = \max(0, \sin(2\pi t/96 - \pi/2))$.  We solve the first of the RP method with this irradiance, and then either lower or raise it it 20\% before running the second stage.  Fig.~\ref{fig:uncer_comp} shows the resulting costs and charges relative to the first stage result.


Lowering solar output produces increases to generation costs of between 8-10\% relative to the original solution.  Function 0 performs around 1\% better on average and it increases the charge to the house in line with the increased generation costs.  For function 3 the households are barely penalised for this increase in generation costs.  When solar output is raised function 3 doesn't react at all because it has no incentive to.  Function 0 on the other hand takes advantage of the excess solar to lower overall costs by around 8\%.

As we found in the previous section, generator costs are relatively independent of the RP cost function penalty parameter $\alpha$.  However, $\alpha$ can be used to increase household charges.  This means that $\alpha$ can be tuned over time to ensure that the market is budget balanced.  For example, these extra penalties can be used to recover frequency regulation costs or to help pay for fixed network costs.

Larger penalties could also help to deter agents looking to game the system.  In the current mechanism, there is nothing that prevents an agent from misreporting their behaviour during the negotiation phase.  A larger penalty forces them to act on what they negotiated, so can significantly limit this behaviour.  The problem is that it will be difficult to distinguish between someone looking to cheat the system, and someone with a large amount of local uncertainty or discrete loads.  This makes it difficult to penalise one of these effects without penalising the others.  Houses with batteries have more flexibility in how they can deal with local sources of uncertainty, so a system with batteries and reasonable penalties might provide the best compromise.

\iftoggle{future}{}{%
We plan to investigate the ability of agents to game the system in future work.
}

\section{Conclusion \iftoggle{future}{and Future Work}{}}
We have improved on existing distributed ADMM based D-OPF methods by including more accurate line models and a two-stage approach to manage discrete variables and uncertainty. We developed a suburb sized test microgrid, and found that the full non-convex AC equations produce close to optimal solutions in short solve times.  Our two-stage approach provides a simple but effective means of managing household discrete variables and uncertainty in the network.

\iftoggle{future}{%
Future research will focus on investigating alternative distributed solving techniques with the aim of further improving convergence.  There are also opportunities to parallelise the problem more by decomposing some components across time.  Further research needs to be done on the real-time control that takes place within and between time-steps, and it might be possible to build a frequency regulation market into our distributed algorithm.

We need further experiments to investigate if our results carry over to larger discrete decisions, for example, those related to large industrial plant, generator start-up costs, and line switching.  We also plan to answer the important question of how susceptible this mechanism is to gaming in practice, and if this is a problem, what can be done about it.
}{}

\Omit{
\appendices
\section*{Acknowledgment}
Thanks to Hassan Hijazi, Matt Kraning.
}

\bibliographystyle{splncs}
\bibliography{../bib/bib}

\end{document}